\documentclass[11pt,reqno]{amsart}
\usepackage{amsthm}
\usepackage{amssymb}
\usepackage{latexsym}
\usepackage{multicol}
\usepackage{verbatim,enumerate}
\usepackage{amscd}

 \usepackage[usenames]{color}
\usepackage{hyperref}
\usepackage{amsmath, amscd}

\advance\textwidth by 1.2in \advance\oddsidemargin by -.6in \advance\evensidemargin by -.6in
\theoremstyle{definition}
\newtheorem{defn}{Definition}[section]
\newtheorem{thm}{Theorem}[section]

\newtheorem{rem}{Remark}[section]

\newtheorem{exam}{Example}[section]
\newtheorem{cor}{Corollary}[section]
\newtheorem{lem}{Lemma}[section]
\newtheorem{prop}{Proposition}[section]

\newcounter{cnt}
 \makeatletter
\def\mydggeometry{\makeatletter\dg@YGRID=1\dg@XGRID=20\unitlength=0.003pt\makeatother}
\makeatother \theoremstyle{remark}

\numberwithin{equation}{section}

\let\bwdg\bigwedge
\def\bigwedge{{\textstyle\bwdg}}

\newcommand{\id}{\operatorname{id}}

\newcommand{\nc}{\newcommand}
\newcommand{\rnc}{\renewcommand}

\nc{\cal}{\mathcal} \nc{\goth}{\mathfrak} \rnc{\bold}{\mathbf}

\renewcommand{\Bbb}{\mathbb}
\nc\bomega{{\mbox{\boldmath $\omega$}}} \nc\bpsi{{\mbox{\boldmath $\Psi$}}}
 \nc\balpha{{\mbox{\boldmath $\alpha$}}}
 \nc\bpi{{\mbox{\boldmath $\pi$}}}
\nc\bmu{{\mbox{\boldmath $\mu$}}} \nc\bcN{{\mbox{\boldmath $\cal{N}$}}} \nc\bcm{{\mbox{\boldmath $\cal{M}$}}} \nc\blambda{{\mbox{\boldmath
$\lambda$}}}\nc\bnu{{\mbox{\boldmath $\nu$}}}

\newcommand{\lie}[1]{\mathfrak{#1}}

\makeatletter
\def\section{\def\@secnumfont{\mdseries}\@startsection{section}{1}%
  \z@{.7\linespacing\@plus\linespacing}{.5\linespacing}%
  {\normalfont\scshape\centering}}
\def\subsection{\def\@secnumfont{\bfseries}\@startsection{subsection}{2}%
  {\parindent}{.5\linespacing\@plus.7\linespacing}{-.5em}%
  {\normalfont\bfseries}}
\makeatother

 \nc{\Hom}{\operatorname{Hom}}
  \nc{\mode}{\operatorname{mod}}
\nc{\End}{\operatorname{End}} \nc{\wh}[1]{\widehat{#1}} \nc{\Ext}{\operatorname{Ext}} \nc{\ch}{\text{ch}} \nc{\ev}{\operatorname{ev}}
\nc{\Ob}{\operatorname{Ob}} \nc{\soc}{\operatorname{soc}} \nc{\rad}{\operatorname{rad}} \nc{\head}{\operatorname{head}}

 \nc{\Cal}{\cal} \nc{\Xp}[1]{X^+(#1)} \nc{\Xm}[1]{X^-(#1)}
\nc{\on}{\operatorname} \nc{\Z}{{\bold Z}} \nc{\J}{{\cal J}} \nc{\C}{{\bold C}} \nc{\Q}{{\bold Q}}

\nc{\N}{{\Bbb N}} \nc\boa{\bold a} \nc\bob{\bold b} \nc\boc{\bold c} \nc\bod{\bold d} \nc\boe{\bold e} \nc\bof{\bold f} \nc\bog{\bold g}
\nc\boh{\bold h} \nc\boi{\bold i} \nc\boj{\bold j} \nc\bok{\bold k} \nc\bol{\bold l} \nc\bom{\bold m} \nc\bon{\bold n} \nc\boo{\bold o}
\nc\bop{\bold p} \nc\boq{\bold q} \nc\bor{\bold r} \nc\bos{\bold s} \nc\boT{\bold t} \nc\boF{\bold F} \nc\bou{\bold u} \nc\bov{\bold v}
\nc\bow{\bold w} \nc\boz{\bold z} \nc\boy{\bold y} \nc\ba{\bold A} \nc\bb{\bold B} \nc\bc{\bold C} \nc\bd{\bold D} \nc\be{\bold E} \nc\bg{\bold
G} \nc\bh{\bold H} \nc\bi{\bold I} \nc\bj{\bold J} \nc\bk{\bold K} \nc\bl{\bold L} \nc\bm{\bold M} \nc\bn{\bold N} \nc\bo{\bold O} \nc\bp{\bold
P} \nc\bq{\bold Q} \nc\br{\bold R} \nc\bs{\bold S} \nc\bt{\bold T} \nc\bu{\bold U} \nc\bv{\bold V} \nc\bw{\bold W} \nc\bz{\mathbb{ Z}} \nc\bx{\bold
x} \nc\KR{\bold{KR}} \nc\rk{\bold{rk}} \nc\het{\text{ht }}

\nc\toa{\tilde a} \nc\tob{\tilde b} \nc\toc{\tilde c} \nc\tod{\tilde d} \nc\toe{\tilde e} \nc\tof{\tilde f} \nc\tog{\tilde g} \nc\toh{\tilde h}
\nc\toi{\tilde i} \nc\toj{\tilde j} \nc\tok{\tilde k} \nc\tol{\tilde l} \nc\tom{\tilde m} \nc\ton{\tilde n} \nc\too{\tilde o} \nc\toq{\tilde q}
\nc\tor{\tilde r} \nc\tos{\tilde s} \nc\toT{\tilde t} \nc\tou{\tilde u} \nc\tov{\tilde v} \nc\tow{\tilde w} \nc\toz{\tilde z} \nc\woi{w_{\omega_i}}
\nc\chara{\operatorname{Char}}
\nc{\buf}{\underline{\bof}} 
\nc{\avee}{\alpha^{\vee}}
\nc{\bLambda}{\mathbf{\Lambda}}
\nc{\bMu}{\mathbf{M}}
\nc{\bNu}{\mathbf{N}}
\nc{\btau}{\mathbf{\sigma}}
\nc{\bsigma}{\mathbf{\sigma}}

\begin{document}

\author[Fourier]{Ghislain Fourier}
\address{Mathematisches Institut, Universit\"at zu K\"oln, Germany}
\email{gfourier@math.uni-koeln.de}
\address{School of Mathematics and Statistics, University of Glasgow, UK}
\email{ghislain.fourier@glasgow.ac.uk}

\thanks{The author was partially supported by the DFG priority program 1388 ``Representation Theory``}

\date{\today}

\title{Extended partial order and applications to tensor products}
\begin{abstract} 
We extend the preorder on $k$-tuples of dominant weights of a simple complex Lie algebra $\lie g$ of classical type adding up to a fixed weight $\lambda$ defined by V.~Chari, D.~Sagaki and the author (\cite{CFS2012}). We show that the induced extended partial order on the equivalence classes has a unique minimal and a unique maximal element. For $k=2$ we compute its size and determine the cover relation.\\ 
To each $k$-tuple we associate a tensor product of simple $\lie g$-modules and we show that for $k=2$ the dimension increases also along with the extended partial order, generalizing a theorem proved in the aforementioned paper. We also show that the tensor product associated to the maximal element has the biggest dimension among all tuples for arbitrary $k$, indicating that this might be a symplectic (resp. orthogonal) analogue of the row shuffle defined by Fomin et al  \cite{FFLP2005}.\\
The extension of the partial order reduces the number elements in the cover relation and may facilitate the proof of an analogue of Schur positivity along the partial order for symplectic and orthogonal types.
\end{abstract}
\maketitle

\section{Introduction}
Let $\lie g$ be a finite-dimensional, simple complex Lie algebra and $P^+$ the set of dominant integral weights. For a given $\lambda \in P^+, k \geq 1$ let $\mathcal{P}^+(\lambda,k)$ denote the subset of $k$-tuples of dominant integral weights adding up to $\lambda$. Chari, Sagaki and the author (\cite{CFS2012}) have defined a preorder $\preceq$ on $\mathcal{P}^+(\lambda,k)$, by extending a natural preorder for $\lie {sl}_2$ and $k=2$, that is if $\lambda = n \omega$, then 
\[ (\lambda - a \omega, a \omega) \preceq (\lambda - b \omega, b \omega) :\Leftrightarrow \min \{n-a, a\} \leq \min \{ n-b, b\}.\] 
To each element $(\lambda_1, \ldots, \lambda_k) \in \mathcal{P}^+(\lambda,k)$, the tensor product of simple  finite-dimensional modules $V(\lambda_1) \otimes \ldots \otimes V(\lambda_k)$ is associated. It was shown that the dimension of this tensor product increases along with the preorder (\cite[Theorem 1(i)]{CFS2012}). Even more, for $\lie g$ of type $A_2$ and $k=2$ or $\lie g$ and $ k$ arbitrary and $\lambda$ a multiple of a minuscule weight, it was shown that there exist injective maps of $\lie g$-modules along with the preorder (\cite[Theorem 1 (ii),(iii)]{CFS2012}).\\
We shall recall that the preorder on $\mathcal{P}^+(\lambda,k)$ depends on the root system of $\lie g$, in fact each positive root gives a certain set of inequalities that determinede the preorder. This allows a quite elementary proof of the dimension theorem using Weyl's dimension formula. On the other hand, to understand the partial order, one is interested in the cover relation. For a given element $(\lambda_1, \ldots, \lambda_k)$ the cover elements are quite hard to compute and it involves the combinatorics of the root system and the associated Weyl group of $\lie g$. In \cite{CFS2012} a classification of covers is given for $k=2$ only.\\
In order to prove the existence of an injective map of $\lie g$-modules along with the preorder,  it is enough to prove it for the cover relation only. The existence was proven for $\lie g$ of type $A_2$ for each possible cover by using the combinatorics of Young tableaux. But even in type $C_2$, the number of possible covers increases as the order of the Weyl group increases, this makes the case by case consideration (as for $A_2$) not appropriate.\\
In the present paper, we extend the partial order and prove the tensor product theorem even for this extended order. For this let $\mathcal{P}^+ = \mathbb{Z}_{\geq 0}^n$ and, as above, for a given $\lambda \in \mathcal{P}^+$, we denote the $k$-tuples adding up to $\lambda$ by $\mathcal{P}^+(\lambda,k)$. We introduce a preorder $\trianglelefteq$ on $\mathcal{P}^+(\lambda,k)$ and we show that if one uses the canonical embedding $\iota:\mathcal{P}^+ \hookrightarrow P^+$, the dominant weights of a simple Lie algebra, then $\trianglelefteq$ is an extension of $\preceq$. In fact if $\lie g$ is of type $A_n$ the partial orders are equivalent (Remark~\ref{an-rem}), we will use this to adapt certain results from \cite{CFS2012} to the present paper.\\

\noindent
We show that $\trianglelefteq$ has a unique maximal element $\blambda^{\max} \in \mathcal{P}^+(\lambda,k)/\!\sim$, the equivalence classes induced by the preorder.  For the preorder $\preceq$ this is true in type $A_n$ only. See for example (Example~\ref{exam-max}) where in type $C_2$ for a certain $\lambda$, the set $\mathcal{P}^+(\iota(\lambda),2)$ has two maximal elements, while $\mathcal{P}^+(\lambda,2)$ has a unique one. Further, we can compute the size of the poset and give a classification of the cover relation.

\noindent
To $(\lambda_1, \ldots, \lambda_k) \in \mathcal{P}^+(\lambda,k)$ we associate a tensor product of simple finite-dimensional $\lie g$-modules, namely $V(\iota(\lambda_1)) \otimes \ldots \otimes V(\iota(\lambda_k))$. We will show
\begin{thm}\label{mainthm}
Let $\blambda = (\lambda_1, \lambda_2) \trianglelefteq \bmu = (\mu_1, \mu_2) \in \mathcal{P}^+(\lambda,2)$, then 
$$\dim V(\iota(\lambda_1)) \otimes V(\iota(\lambda_2)) \leq \dim V(\iota(\mu_1)) \otimes V(\iota(\mu_2)),$$
with equality iff $\blambda = \bmu \in  \mathcal{P}^+(\lambda,2)/\!\sim$.
\end{thm}

\noindent
This theorem was proven in \cite{CFS2012} with respect to $\preceq$, by using Weyl's dimension formula. It was enough to prove it for the $\lie{sl}_2$-case where this is quite elementary. In our case, the preorder $\trianglelefteq$, we have to prove the case where $\lie g$ is of type $C_2$ separately using a case by case consideration (Theorem~\ref{C2}). From here and the $\lie{sl}_2$-case we can deduce the proof in general (Section~\ref{section-proof}).

\noindent
Using the theorem we can show that for all $\bmu = (\mu_1, \ldots, \mu_k) \in \mathcal{P}^+(\lambda, k)$  we have
\begin{eqnarray}
\dim V(\iota(\mu_1)) \otimes \ldots \otimes V(\iota(\mu_k)) \leq \dim V(\iota(\lambda_1)) \otimes \ldots \otimes V(\iota(\lambda_k)) 
\label{eq:max-dim}
\end{eqnarray}
where $ (\lambda_1,\ldots, \lambda_k) = \blambda^{\max}$ and equality iff $\bmu = \blambda$ in $\mathcal{P}^+(\lambda,k)/\!\sim$.\\

\noindent
We shall make a couple of remarks here. If $\lie g $ is of type $A_n$, then \eqref{eq:max-dim} follows from the stronger statement that there is an injective map of $\lie g$-module from the tensor product of left hand side of (\ref{eq:max-dim}) to the tensor product of the right hand side. This statement was proven via the strong connection from $\blambda^{\max}$ to the so-called row shuffle (see \cite[Section 2.4, 2.5]{CFS2012} for more details). The existence of such an injective map was conjectured by Fomin-Fulton-Li-Poon, Lascoux-Leclrec-Thibon and Okounkov (\cite{FFLP2005}, \cite{O1997}, \cite{LLT1997}), and recently proven by Lam, Postnikov, Pylyavskyy (\cite{LPP2007}). A phenomenon like this is called Schur positivity (the difference of the characters of the modules is a non-negative linear combination of Schur functions). For more on this subject see also \cite{BBR2006}, \cite{BM2004}, \cite{DP2007}, \cite{Mc2008}, \cite{MW2009}, \cite{PW2008}. One might see (\ref{eq:max-dim}) as evidence that $\blambda^{\max}$ could be the generalization of the row shuffle to simple Lie algebras of symplectic or orthogonal type.\\
We conjecture that we have injective maps of $\lie g$-modules along with $\trianglelefteq$. Since the order $\trianglelefteq$ reduces the number of possible covers enormously, we have reduced the cases that have to be considered in order to prove a kind of Schur positivity in the symplectic or orthogonal case. This will be discussed elsewhere.

\noindent
The paper is organized as follows:\\
In Section~\ref{section-def} we recall necessary notations of simple Lie algebras, define the preorder $\trianglelefteq$ and relate it to the preorder $\preceq$. In Section~\ref{section-pre}  we analyze the preorder further, computing the unique maximal element, while in Section~\ref{section-pre2} we restrict ourselves to $k=2$ and recall the previously known results about the cover relation. Section~\ref{section-proof} is dedicated to the proof of Theorem~\ref{mainthm}.


\section{Notations and definitions}\label{section-def}
In this section we recall the basic notations and introduce the main object of this paper, the preorder $\trianglelefteq$.\\

\noindent
Let $\mathcal{P} = \bz^n$ and $\mathcal{P}^+ = \bz_{\geq 0}^n$, we set $I = \{1, \ldots, n\}$. We denote $\{ \omega_i = (0,\ldots,0,1,0, \ldots,0) \, | \, i \in I \}$ the canonical basis of $\mathcal{P}$. We denote further $\omega_i^* \in \Hom_{\bz}(\mathcal{P}, \bz)$ the dual element of $\omega_i$. Then for any $\lambda \in \mathcal{P}$ we have 
$$\lambda = \sum_{i \in I} \omega_i^*(\lambda) \omega_i.$$ 
We have an alternative description by setting $\epsilon_i = \omega_i - \omega_{i-1}$, then 
$$\lambda = \sum_{i \in I} (\omega_i^* + \ldots + \omega_n^*)(\lambda) \epsilon_i.$$\\
Let $S_n$ be the symmetric group in $n$-letters and $s_{i, i+1}$ ($ i =1 , \ldots , n-1$) the simple transpositions generating $S_n$. $S_n$ acts on $\mathcal{P}$ via
$$s_{i,i+1} (\epsilon_j) = \epsilon_{s_{i,i+1}(j)} = \begin{cases} \epsilon_j \text{ if } j < i \text{ or } j > i+1 \\ \epsilon_{j+1} \text{ if } j = i \\ \epsilon_{j-1} \text{ if } j = i+1 \end{cases}
$$
It is standard to see that for any $\lambda = (a_1,  \ldots, a_n) \in \mathcal{P}$
\[ |S_n \lambda \cap \mathcal{P}^+| = 1 ,\]
so there is a unique element in the $S_n \lambda$-orbit which is in $\mathcal{P}^+$.\\

\subsection{}
We fix $\lambda \in \mathcal{P}^+$ and $k \in \bz_{\geq 0}$, and we set 
$$\mathcal{P}^+(\lambda,k) = \{ \blambda = (\lambda_1, \ldots, \lambda_k) \, | \, \lambda_i \in \mathcal{P}^+ \, , \, \lambda_1 + \ldots + \lambda_k = \lambda \}.$$

Then $\mathcal{P}^+(\lambda,k)$ is a finite set for all $\lambda \in \mathcal{P}^+$. Given $\blambda \in \mathcal{P}^+(\lambda,k)$ we define for each pair $i \leq j \in I$ and $1 \leq  \ell \leq k$: 
$$\bor_{(i,j), \ell}(\blambda) = \min \{ \sum_{ k = i}^{j} \omega_k^*(\lambda_{n_1}  + \ldots + \lambda_{n_\ell}) \, | \, 1 \leq n_1 < \ldots < n_\ell \leq k \}.$$
We have for all $i \leq j \in I,\ell,\blambda$
$$0 \leq \bor_{(i,j),\ell}(\blambda) \leq \bor_{(i,j), \ell + 1}(\blambda),$$
and for all $i \leq j \in I$ and $\blambda, \bmu \in \mathcal{P}^+(\lambda,k)$
\[ \bor_{(i,j),k}(\blambda) = \bor_{(i,j), k}(\bmu).\]
Given $\blambda, \bmu \in \mathcal{P}^+(\lambda,k)$, we set 
\begin{eqnarray}
\blambda \trianglelefteq \bmu :\Leftrightarrow \bor_{(i,j), \ell}(\blambda) \leq \bor_{(i,j), \ell}(\bmu)  \text{ for all } i \leq j \in I, 1\leq  \ell \leq k.
\label{eq:ord-def}
\end{eqnarray}
This defines a preorder on $\mathcal{P}^+(\lambda,k)$ and we set
$$\blambda \sim \bmu :\Leftrightarrow \bor_{(i,j), \ell}(\blambda) = \bor_{(i,j), \ell}(\bmu)  \text{ for all }  i \leq j \in I, 1\leq  \ell \leq k,$$
to obtain an induced partial order on $\mathcal{P}^+(\lambda,k)/\!\sim$.

\noindent
$S_k$ acts on $\mathcal{P}^+(\lambda,k)$ by permuting the components but is invariant on equivalence classes:
$$ \sigma(\blambda) \sim \blambda \; \forall \; \sigma \in S_k.$$

\subsection{}
Following Bourbaki \cite{Bou75}, we introduce certain notations for Lie algebras.

\noindent
Let $\lie g$ be a simple, complex finite dimensional Lie algebra of classical type, denote $n$ the rank of $\lie g$, $I =  \{1, \ldots, n\}$ and fix a triangular decomposition $\lie n^+ \oplus \lie h \oplus \lie n^-$. Denote $R$ (resp. $R^+$) the (positive roots) roots of $\lie g$, $P$ (resp. $P^+$) the (dominant) integral weights. For a given $\alpha \in R^+$, the coroot is denoted $h_\alpha$. The simple roots are denoted $\alpha_i$, the corresponding coroots $h_i$ and the fundamental dominant weights $\omega_i$. We denote $W$ the Weyl group associated to $\lie g$, $\bos_i$ the simple reflections and $(\,, \,)$ an invariant bilinear form on $\lie h^* \times \lie h$. For later use we will give a complete list of all positive coroots (in terms of simple coroots).
\begin{rem}\label{coroot-list} List of coroots of simple, finite-dimensional Lie algebras of classical type.
\begin{enumerate}
\item type $A_n$: \\ $$\{ h_{i,j} = h_{i} + \ldots + h_j \, | \,  i \leq j \in I \}$$
\item type $C_n$:\\ $$\{ h_{i,j} = h_{i} + \ldots + h_j \, | \,  i \leq j \in I \}$$  
$$ \{ h_{i,\overline{j}} = h_{i} + \ldots + h_{j-1} + 2 h_{j} + \ldots + 2 h_{n} \, | \, i < j \in I \}$$
\item type $B_n$:\\ $$\{ h_{i,j} = h_{i} + \ldots + h_j \, | \, i \leq j \in I \}$$  
$$ \{ h_{i,\overline{j}} = h_{i} + \ldots + h_{j-1} + 2 h_{j} + \ldots + 2 h_{n-1} + h_{n} \, | \,  i \leq j \in I\}$$
\item type $D_n$:\\ $$\{ h_{i,j} = h_{i} + \ldots + h_j \, | \,  i \leq j \in I \}$$ 
$$  \{ h_{i, \overline{n}} = h_{i} + \ldots + h_{n-2} + h_{n} \, | \, 1 \leq i \leq n-2 \}$$ 
$$ \{ h_{i,\overline{j}} = h_{i} + \ldots + h_{j-1} + 2 h_{j} + \ldots + 2 h_{n-2}  +  h_{n-1} + h_{n} \, | \,  1 \leq i < j < n-1 \}.$$
\end{enumerate}
\end{rem}


\section{About the partial order}\label{section-pre}
In \cite[Section 2.1]{CFS2012} a preorder $\preceq$ on $k$-tuples of dominant weights for finite-dimensional, simple complex Lie algebra $\lie g$ of finite rank was introduced. We recall this definition here. For a given $\lambda \in P^+$, $k \in \bz_{\geq 1}$ we set
$$P^+(\lambda,k) = \{ \blambda = (\lambda_1, \ldots, \lambda_k) \, | \, \lambda_i \in P^+ \, , \, \lambda_1 + \ldots + \lambda_k = \lambda \}.$$
Let $\blambda \in P^+(\lambda,k)$, the following integer was defined for $1 \leq \ell \leq k$ and  $\alpha \in R^+$
$$\bor_{\alpha, \ell}(\blambda) = \min \{ (\lambda_{i_1} + \ldots + \lambda_{i_{\ell}})(h_\alpha) \, | \, 1 \leq i_1 < \ldots < i_{\ell} \leq k \}.$$
The preorder $\preceq$ was defined as follows:\\
Let $\blambda, \bmu \in P^+(\lambda,k)$, then 
\begin{eqnarray}
\blambda \preceq \bmu :\Leftrightarrow  \bor_{\alpha, \ell}(\blambda) \leq \bor_{\alpha, \ell}(\bmu)  \text{ for all } \alpha \in R^+, 1\leq  \ell \leq k.
\label{eq:prec-def}
\end{eqnarray}

\noindent
We want to identify $\mathcal{P}^+$ with a certain subset of $P^+$. Let $\lie g$ be of type $A_n, C_n, B_{n+1}, D_{n+2}$, then we embed $\iota: \mathcal{P}^+ \hookrightarrow P^+$ by $\omega_i \mapsto  \omega_i$ for $1 \leq i \leq n$. 
We call the weights in $\iota(\mathcal{P}^+)$ admissible, e.g. every $\lambda \in P^+$ is admissible if $\lie g $ is of type $A_n, C_n$, $\lambda$ is admissible if $\lambda(h_{n+1}) = 0$ if $\lie g$ is of type $B_{n+1}$ or $\lambda(h_{n+1}), \lambda(h_{n+2}) = 0$ if $\lie g$ is of type $D_{n+2}$. In other words, a weight is called admissible if it is not supported on a spin node. We have an induced embedding $\iota: \mathcal{P}^+(\lambda, k) \hookrightarrow P^+(\lambda,k)$.\\
Let $\lambda \in \mathcal{P}^+$ and $\alpha \in R^+$, such that $h_\alpha = h_i + \ldots + h_j$ for some $i \leq j$, then 
\[ \bor_{\alpha, \ell}(\iota(\blambda)) = \bor_{(i,j), \ell} (\blambda).\]
This implies immediately
\begin{cor}\label{cor-compare}
Let $\lie g$ be a complex simple Lie algebra, $\lambda \in \mathcal{P}^+$ and $\blambda, \bmu \in \mathcal{P}^+(\lambda,k)$
\[ \iota(\blambda) \preceq \iota(\bmu) \Rightarrow \blambda \trianglelefteq \bmu.\]
\end{cor}

\noindent
\begin{rem}\label{an-rem}
If $\lie g$ is of type $A_n$ all coroots are of the form $h_i+ \ldots + h_j$ and, since $\iota$ is surjective, all weights are admissible. In this case we have
\[ \iota(\blambda) \preceq \iota(\bmu) \Leftrightarrow \blambda \trianglelefteq \bmu.\]
We will use this fact throughout the article to adapt certain results of \cite{CFS2012}.
\end{rem}

\subsection{}
Given $\lambda \in \mathcal{P}^+$, then it is easy to see that the poset $\mathcal{P}^+(\lambda,k)/\!\sim$ has a unique minimal element $\blambda^{\min}$, namely the $S_k$-orbit of $\blambda = (\lambda, 0, \ldots, 0)$. We will show that there is also a unique maximal element $\blambda^{\max}$.\\
Let $\lambda = (a_1, \ldots, a_n) \in \mathcal{P}^+$, then for $i \in I$
\[ \sum_{\ell =i}^n a_{\ell} = p_i k + r_i \; \text{ for some }0 \leq p_i, 0 \leq r_i < k.\]
We set for $i \in I , 1 \leq j \leq k$
$$ m^{i,j} = 
\begin{cases}
p_i + 1 & \text{ if } j \leq r_i \\ 
p_i & \text{ if } j > r_i
\end{cases}
$$
then $m^{i,j} \geq m^{i+1,j}$ for any $i,j$. This implies that $\lambda_j := \sum\limits_{i = 1}^n m^{i,j} \epsilon_i \in \mathcal{P}^+$.
We set 
\begin{eqnarray}
\blambda^{\max} := (\lambda_1, \ldots, \lambda_k).
\label{eq:max-def}
\end{eqnarray}
\begin{lem}\label{maximal-element}
Let $\bmu \in \mathcal{P}^+(\lambda,k)$, then $\bmu \trianglelefteq \blambda^{\max}$. And if $\blambda^{\max} \sim \bmu$ then $\bmu$ and $\blambda^{\max}$ are in the same $S_k$-orbit.
\end{lem}

Before proving the lemma, we will show a useful proposition. 
\begin{prop}\label{prop:extend-tuple} Let $ \mu, \tau \in \mathcal{P}^+$, $s < k \in \mathbb{Z}_{\geq 1}$. Suppose $(\lambda_1, \ldots, \lambda_s) \triangleleft (\mu_1, \ldots, \mu_s) \in \mathcal{P}^+(\mu, s)$ and $(\tau_1, \ldots, \tau_{k-s}) \in \mathcal{P}^+(\tau, k-s)$, then
\[
(\lambda_1, \ldots, \lambda_s, \tau_1, \ldots, \tau_{k-s}) \triangleleft (\mu_1, \ldots, \mu_s, \tau_1, \ldots, \tau_{k-s}) \in \mathcal{P}^+(\mu + \tau, k).
\]
\end{prop}
\begin{proof}
First we reduce the proof to the case $\mathcal{P}^+ = \mathbb{Z}_{\geq 0}$ (e.g. $n =1$). For this reduction we define for each $i \leq j$ a map
$$\pi_{i,j} : \mathcal{P}^+ \longrightarrow \mathbb{Z}_{\geq  0} \; : \; (a_1, \ldots, a_n) \mapsto a_i + \ldots + a_j.$$
We can extend this component wise to a map 
\begin{eqnarray}
\pi_{i,j} :\mathcal{P}^+ (\lambda, k) \longrightarrow \mathbb{Z}_{\geq 0}(\pi_{i,j}(\lambda), k).
\label{eq:pi-ij}
\end{eqnarray}
Let $\lambda \in \mathcal{P}^+$ and  $\blambda \trianglelefteq \bmu \in \mathcal{P}^+(\lambda, k)$, then 
$$\pi_{i,j} (\blambda) \trianglelefteq \pi_{i,j}(\bmu) \in  \mathbb{Z}_{\geq 0}(\pi_{i,j}(\lambda), k).$$
Even more if $\blambda \triangleleft \bmu$, then there exists $i \leq j$ such that $\pi_{i,j} (\blambda) \triangleleft \pi_{i,j}(\bmu) \in  \mathbb{Z}_{\geq 0}(\pi_{i,j}(\lambda), k)$.\\
On the other hand if $\pi_{i,j} (\blambda) \trianglelefteq \pi_{i,j}(\bmu) \in  \mathbb{Z}_{\geq 0}(\pi_{i,j}(\lambda), k)$ for all $i \leq j$, then $\blambda \trianglelefteq \bmu$, with equality iff there is equality for all $i \leq j$. So we have (a reformulation of \ref{eq:ord-def}):
\begin{eqnarray}
\blambda \trianglelefteq \bmu \Leftrightarrow \pi_{i,j}(\blambda) \trianglelefteq \pi_{i,j}(\bmu) \text{ for all } i \leq j
\label{eq:pi-equiv}
\end{eqnarray}
which reduces the proof to the case $n=1$. Let $\mathcal{P}^+ = \mathbb{Z}_{\geq 0}$, $s < k$, $\mu, \tau \in \mathcal{P}^+$, $(a_1, \ldots, a_s) \triangleleft (b_1, \ldots, b_s) \in \mathcal{P}^+(\mu, s)$, where $a_i, b_i \in \mathbb{Z}_{\geq 0}$. We may assume $a_1 \leq \ldots \leq a_s$ and $b_1 \leq \ldots \leq b_s$, then by assumption 
\begin{eqnarray}
\sum_{j = 1}^{i} a_j \leq \sum_{j = 1}^i b_j \text{ for all } 1 \leq i \leq s
\label{eq:3-1-a}
\end{eqnarray}
and there exists $1 \leq i < s$ such that the inequality is strict.\\
Let $c \in  \mathbb{Z}_{\geq 0}$ and $i_1 \in \{1, \ldots, s\}$ be the minimum such that $c < a_{i_1}$, $i_2$ the minimum such that $c < b_{i_2}$. We want to check that for all $1 \leq i \leq s+1$
\begin{eqnarray}
r_{(1,1), i} (a_1, \ldots, a_s, c) \leq r_{(1,1), i} (b_1, \ldots, b_s, c).
\label{eq:max-1}
\end{eqnarray}
and there exists $1 \leq i < s+1$ such that the inequality is strict.\\
For $i \leq \min\{i_1, i_2 \}$ or $i \geq \max\{i_1, i_2\}$ this follows by assumption. To complete this, we must distinguish two more cases, $i_1 \leq i \leq i_2$ and $i_2 \leq i \leq i_1$. About the first one, the left hand side of \ref{eq:max-1} is 
$$c + \sum_{j=1}^{i -1} a_j < \sum_{j=1}^{i} a_j \leq \sum_{j=1}^{i} b_j$$
(the first inequality is due to the definition of $i_1$ and the second due to the assumptions), which is the right hand side of the equation.\\
Let now $i_2 \leq i \leq i_1$, then the left hand side of \ref{eq:max-1} is
$$\sum_{j=1}^{i} a_j \leq c + \sum_{j=1}^{i-1} a_j \leq c + \sum_{j=1}^{i-1} b_j <  \sum_{j=1}^{i} b_j,$$
(where the first inequality is due to the definition of $i_1$, the second by assumption and the third due to the definition of $i_2$), which is the right hand side of \ref{eq:max-1}. If $\ell$ is such that the inequality in \eqref{eq:3-1-a} is strict then the inequality is strict for $\ell -1, \ell$ or $\ell +1$ resp.\\
So we have proven that
$$(a_1, \ldots, a_s) \triangleleft (b_1, \ldots, b_s) \Rightarrow (c, a_1, \ldots, a_s) \triangleleft (c, b_1, \ldots, b_s).$$
By induction we see
$$(a_1, \ldots, a_s) \triangleleft (b_1, \ldots, b_s) \Rightarrow (c_1, \ldots, c_{k-s}, a_1, \ldots, a_s) \triangleleft (c_1, \ldots, c_{k-s}, b_1, \ldots, b_s).$$
\end{proof}

\begin{proof}(of Lemma~\ref{maximal-element})
The lemma has been already proven for $\mathcal{P}^+= \mathbb{Z}_{\geq 0}^n$ for arbitrary $n$ and $k=2$ \cite[Proposition 5.3]{CFS2012} as well as for $n=1$ and arbitrary $k$ \cite[Lemma 3.3]{CFS2012}.  It was shown there that if $(a_1, \ldots, a_k)$ is a maximal element in $\mathbb{Z}_{\geq 0}(\lambda,k)$, then $a_i - a_j \in \{0, \pm 1\}$ for all $i \leq j \in I$. By ordering the tuples (recall the $S_k$-action), the unique ordered maximal element is determined by the condition
 \begin{eqnarray}
a_i - a_j \in \{0,1\} \text{ for }i \leq j \in I
\label{eq:max-uni}
\end{eqnarray}
Let $\blambda^{\max} \in \mathcal{P}^+(\lambda,k)$ as defined in \ref{eq:max-def}. Then $\pi_{i,j}(\blambda^{\max}) = (a_1, \ldots, a_k)$ satisfies (\ref{eq:max-uni}) for all $i \leq j \in I$. This implies that for given $\bmu \in \mathcal{P}^+(\lambda, k)$
$$\pi_{i,j}(\bmu) \trianglelefteq \pi_{i,j}(\blambda^{\max}) \text{ for all  } i \leq j \in I.$$
With (\ref{eq:pi-equiv}) we can conclude that $\bmu \trianglelefteq \blambda^{\max}$. It remains to show that if $\blambda^{\max} \sim \bmu$, then $\blambda^{\max}, \bmu$ are in the same $S_k$-orbit.\\
For this we recall that any $\lambda \in \mathcal{P}^+$ can be uniquely written $\lambda =\sum_{i \in I} b_i \epsilon_i$, where $ b_i \geq b_{i+1} \geq 0$ for all $i$. Given $\bnu = ( \nu_1, \ldots, \nu_k) \in \mathcal{P}^+(\lambda, k)$, we order the $\nu_i$ with respect to the lexicographic order (by writing them in terms of the basis $\{ \epsilon_i \, | \, i \in I\}$).\\ 
We will need that $\blambda^{\max}$ (see \ref{eq:max-def}) is uniquely determined by the condition: 
\begin{eqnarray}
\forall \, i \leq j \; : \; \lambda_i - \lambda_j = \sum_{i \in I} c_\ell^{i,j} \epsilon_\ell \text{ with }c_\ell^{i,j} \in \{0,1\}.
\label{eq:max-uni-2}
\end{eqnarray}
First of all, $\blambda^{\max}$ satisfies this condition. On the other hand, writing $\lambda = \sum_{i \in I} b_i \epsilon_i$ and let $b_i = r_i k + p_i$, then the condition implies that $(\omega_i^* + \ldots + \omega_j^*)(\lambda_\ell) \in \{ r_i, r_i+1\}$, which implies that $\blambda^{\max}$ is uniquely determined.\\

\noindent
Given $\bmu = (\mu_1, \ldots, \mu_k) \in \mathcal{P}^+(\lambda,k)$ a maximal element and we may assume that $\mu_1 \geq \ldots \geq \mu_k$ with respect to the lexicographic order. Then $(\mu_i, \mu_j) \in \mathcal{P}^+(\mu_i + \mu_j, 2)$ and suppose it is not the maximal element, then there exists $(\tau_1, \tau_2) \in \mathcal{P}^+(\mu_i + \mu_j, 2)$ such that 
$(\mu_i, \mu_j) \triangleleft (\tau_1, \tau_2)$. Hence if we replace in $\bmu$, $\mu_i$ by $\tau_1$ and $\mu_j$ by $\tau_2$, and denote the obtained $k$-tuple $\bmu'$, then by Proposition~\ref{prop:extend-tuple}
$$\bmu \triangleleft \bmu' \in \mathcal{P}^+(\lambda, k). $$
This is a contradiction, because $\bmu$ is chosen to be maximal. So we have $(\mu_i, \mu_j)$ is maximal in $\mathcal{P}^+(\mu_i + \mu_j, 2)$ for all $i \leq j$. As mentioned before, the case $k=2$ was proven in  \cite[Proposition 5.3]{CFS2012}. It was shown that $\mathcal{P}^+(\mu_i + \mu_j, 2)$ has a unique (up to $S_2$ action) maximal element. This implies that $(\mu_i, \mu_j)$ satisfies (\ref{eq:max-uni-2}) (using the uniqueness for $\mathcal{P}^+(\mu_i + \mu_j, 2)$). So for all $i \leq j$ we have 
$$\mu_i - \mu_j = \sum_{i \in I} c_\ell^{i,j} \epsilon_\ell \text{ with }c_\ell^{i,j} \in \{0,1\}.$$
But then again with (\ref{eq:max-uni-2}) we have $\bmu = \blambda$.
\end{proof}

We shall remark, that a maximal element in $\mathcal{P}^+(\lambda,k)$ with respect to the partial order $\preceq$ as defined in \cite{CFS2012} is not unique in general.  See for example:
\begin{exam}\label{exam-max} Let $k=2$, $n = 2$ and $\lambda = 2\omega_1 + \omega_2$. Let $\lie g$ be of type $C_2$. Then $P^+(\iota(\lambda),2)$ consists of three $S_2$-orbits and we have 
$$(\iota(\lambda),0 ) \prec (2\omega_1, \omega_2) \; \text{ and } (\iota(\lambda), 0) \prec (\omega_1 + \omega_2, \omega_1)$$
but $(2\omega_1, \omega_2)$ and $(\omega_1 + \omega_2, \omega_1)$ are incomparable with respect to $\preceq$. The partial order $\trianglelefteq$ of $\mathcal{P}^+(\lambda,2)$ gives
$$(\lambda,0) \triangleleft (2\omega_1, \omega_2) \triangleleft (\omega_1 + \omega_2, \omega_1).$$
\end{exam}
The extension of the partial order (see Corollary~\ref{cor-compare}) allows us to order the maximal elements (with respect to $\preceq$) to get a unique maximal element with respect to $\trianglelefteq$.


\section{The partial order for \texorpdfstring{$k=2$}{ k equals 2 }}\label{section-pre2}
Throughout this section we will restrict ourselves to $k=2$. Mainly we are using certain results from \cite{CFS2012}. We will investigate on the partial order further, determined the cover relation and the size of the poset $\mathcal{P}^+(\lambda,2)/\sim$.
\begin{lem}\label{equiv-classes}
Suppose $\blambda = (\lambda_1, \lambda_2) \sim \bmu = (\mu_1, \mu_2)  \in \mathcal{P}^+(\lambda,2)$, then either $\mu_1 = \lambda_1, \mu_2 = \lambda_2$ or $\mu_1 = \lambda_2, \mu_2 = \lambda_1$, hence the equivalence classes in $\mathcal{P}^+(\lambda,2)/\!\sim$ are precisely the $S_2$ orbits.
\end{lem}
\begin{proof}
By using Remark~\ref{an-rem}, the proof is analogue to the proof of \cite[Lemma 5.5]{CFS2012} and $\lie g$ being of type $A_n$.
\end{proof}

We denote the equivalence class of $\blambda = (\lambda_1, \lambda_2)$ again by $\blambda$.

\subsection{}
It is useful in understanding the poset $\mathcal{P}^+(\lambda,2)/\!\sim$, to understand the associated cover relation. 

\begin{defn}
Let $\blambda, \bmu \in \mathcal{P}^+(\lambda,2)/\!\sim$, then $\bmu$ is called a cover of $\blambda$ if and only if
\begin{enumerate}
\item $\blambda \triangleleft \bmu$
\item For any $\bnu \in\mathcal{P}^+(\lambda,2)/\!\sim$ with $\blambda \trianglelefteq \bnu \trianglelefteq \bmu$, either $\blambda = \bnu$ or $\bnu = \bmu$.
\end{enumerate}
\end{defn}

For any $\lambda \in \mathcal{P}^+$, $\mathcal{P}^+(\lambda,2)$ is a finite set, hence if $\blambda \triangleleft \bmu$ there exists a finite chain of successive covers 
$$\blambda \triangleleft \bnu_1 \triangleleft \ldots \triangleleft \bnu_k \triangleleft \bmu.$$

\begin{prop}\label{cover}
Let $\lambda \in \mathcal{P}^+$, $\blambda = (\lambda_1, \lambda_2) $ and $\bmu = (\mu_1, \mu_2) \in \mathcal{P}^+(\lambda,2)/\!\sim$.  Let $\sigma \in  S_n$ such that $\sigma(\lambda_1 - \lambda_2) \in \mathcal{P}^+$. Then
\[\blambda \trianglelefteq \bmu \Leftrightarrow \btau(\lambda_1 - \mu_1), \btau(\mu_1 - \lambda_2) \in \mathcal{P}^+.\]
Further if $\blambda \triangleleft \bmu$ is a cover and there exists $i \in I$ such that 
$$\omega_i^*(\sigma(\lambda_1 - \mu_1)) > 0, \omega_i^*(\btau(\mu_1 - \lambda_2)) > 0,$$ 
then 
\[(\mu_1, \mu_2) =  (\lambda_1 - \btau^{-1} \omega_i, \lambda_2 + \btau^{-1} \omega_i)\]
\end{prop}
\begin{proof} Recall here, that the pair $\mathcal{P}, S_n$ can be identified (via the embedding $\iota:\mathcal{P} \longrightarrow P$) with the lattice of integral weights of the simple Lie algebra $\lie{sl}_{n+1}$ and its Weyl group. This implies that the statement follows immediately from \cite[Proposition 5.4]{CFS2012} and Remark~\ref{an-rem}.
\end{proof}

Suppose a cover $\blambda \triangleleft \bmu$ satisfies the conditions in Proposition~\ref{cover}, that is that  there exists $i \in I$ such that 
$$\omega_i^*(\sigma(\lambda_1 - \mu_1)) > 0, \omega_i^*(\btau(\mu_1 - \lambda_2)) > 0,$$
then we call it a \textit{ cover of type I}.\\ 
Suppose now that for all  $i \in I$ 
\[\omega_i^*(\sigma(\lambda_1 - \mu_1)) = 0 \text{ or }\omega_i^*(\btau(\mu_1 - \lambda_2))= 0,\]  
note that by Proposition~\ref{cover} both are $\geq 0$. Then this implies that
$$\bmu = \sigma^{-1} \left(\sum_{i =1}^n \omega_i^*( \sigma(\lambda_{\epsilon_i})) \omega_i\right),$$
for some $\epsilon = (\epsilon_1, \ldots , \epsilon_n) \in \{1,2\}^{\times n}$. In this case we call $\blambda \triangleleft \bmu$ a \textit{cover of type II}.

\begin{cor}\label{number-cover}
Let $\blambda \in \mathcal{P}^+(\lambda,2)/\!\sim$, then $\blambda$ has at most $n$ covers of type I and $2^{n-1}-1$ covers of type II.
\end{cor}
\begin{proof}
For a given $\blambda$, the element $\sigma \in S_n$ such that $\sigma(\lambda_1 - \lambda_2) \in \mathcal{P}^+$ is uniquely determined up to elements from the stabilizer of $\lambda_1 - \lambda_2$. Then it follows that the type I covers are all obtained via this unique $\sigma$. This implies that there are at most $n$ type I covers.\\
For the type II cover, the element $\sigma$ is again fixed, this gives $2^n$ possible covers, the $S_2$ action of $(\lambda_1, \lambda_2)$ reduces this to $2^{n-1}$. Since the identity is trivial, hence not a cover, we have at most $2^{n-1}-1$ covers of type II.
\end{proof}

\subsection{}
In this section we will compute the size of the poset $\mathcal{P}^+(\lambda,2)/\!\sim$. In this case, $k=2$, this is a simple calculation using Burnside's Lemma. For arbitrary $k$ it is still Burnside's Lemma but the computation is more difficult. There is no formula known to the author computing this for $k > 2$ and it seems to be an interesting combinatorial question.

\begin{prop}
Let $\lambda = \sum_{i \in  I} m_i \omega_i \in P^+$. Then
$$
|\mathcal{P}^+(\lambda,2)/\!\sim| = \begin{cases} \frac{1}{2}\left( \prod_{i \in I} (m_i + 1)  + 1 \right)\text{ if all } m_i \text{ are even } \\ \frac{1}{2}\prod_{i \in I} (m_i + 1) \text{ else } \end{cases}
$$
\end{prop}
\begin{proof}
We will use Burnside's Lemma to compute the orbits of the $S_2$-action. So we have to compute the fixed points for the elements of $S_2$ on the set of ordered tuples. The number of ordered tuples is 
\[
\prod_{i \in I} (m_i + 1).
\]
So the number of fixed points of the identity is $\prod_{i \in I} (m_i + 1)$. The non-trivial element in $S_2$ has a fixed point iff all $m_i$ are even, namely the element $(\sum_{i \in I} m_i/2 \omega_i, \sum_{i \in I} m_i/2 \omega_i)$. Then Burnside's formula gives the proposition.
\end{proof}

\subsection{}
We will induce the proof of Theorem~\ref{mainthm} from the rank $2$ case (that is $n = 2$). In order to prove the rank $2$ case, we need more detailed information about the cover relation in this case. The following proposition is adapted from \cite[Proposition 6.1] {CFS2012} and can be proven similarly (by using the identification $\iota$ in the $\lie{sl}_{n+1}$-case).
\begin{prop} \label{cover-elements} 
Let $\lambda \in \mathcal{P}^+$, $\blambda = (\lambda_1, \lambda_2) \in\mathcal{P}^+(\lambda,2)/\!\sim$. If $\blambda \triangleleft \bmu$ is a cover then $\bmu$ is one of the following list: 
\begin{gather} \label{limitbmu}
\begin{cases}
\bmu = (\lambda_{1}-\btau \omega_{1},\ \lambda_{2}+\btau \omega_{1}) \quad \text{ if } \omega_1^*(\btau(\lambda_1 - \lambda_2)) \geq 2 \text{ or } \\[1.5mm]
\bmu = (\lambda_{1}-\btau \omega_{2},\ \lambda_{2}+\btau \omega_{2}) \quad \text{ if } \omega_2^*(\btau(\lambda_1 - \lambda_2)) \geq 2 \text{ or } \\[1.5mm]
\bmu=(\lambda_1 - \omega_1^*(\btau(\lambda_1-\lambda_2)) \btau \omega_1, \
 \lambda_2 + \omega_1^*(\btau(\lambda_1-\lambda_2))\btau \omega_1) \text{ if } \omega_1^*(\btau(\lambda_1 - \lambda_2)) > 0  
\end{cases}
\end{gather}
where $\sigma \in \{ \id, \bos_{12,}, \bos_{2,3} \}$ such that $\sigma(\lambda_1 - \lambda_2) \in \mathcal{P}^+$.
\end{prop}


\section{Application to representation theory}\label{section-proof}
We will apply the results on the partial order to certain tensor products of simple finite--dimensional modules of a simple complex Lie algebra $\lie g$.

\subsection{}
$P^+$ parameterizes the simple finite-dimensional $\lie g$-modules, denote by $V(\lambda)$ the simple module associated to $\lambda \in P^+$. Its dimension is given by Weyl's dimension formula
\[ \dim V(\lambda) = \prod\limits_{\alpha \in R^+} \frac{(\lambda + \rho, h_\alpha)}{(\rho, h_\alpha)} \]
where $\rho$ is half the sum of all positive roots and $( \, , \, )$ is the invariant bilinear form on $\lie h^* \times \lie h$.
If we denote $(\lambda + \rho, h_\alpha)$ by $\langle \lambda, h_\alpha \rangle$, then we obtain that for $\lambda, \mu \in P^+$
\begin{eqnarray}
\dim V(\lambda) \otimes V(\mu) = \prod\limits_{\alpha \in R^+} \langle \lambda, h_\alpha \rangle \, \langle \mu, h_\alpha \rangle \, (\rho, h_\alpha)^{-2} 
\label{eq:dim-tensor}
\end{eqnarray}

\subsection{}
We want to prove Theorem~\ref{mainthm} for $\mathcal{P}$ of rank $2$, that is $n=2$. For this we need the following useful lemma:
\begin{lem}\label{4-tuple}
Let $a,b,c,d \in \bz_{>0}$, such that $a < b < d$, $a < c < d$, and $b-a \geq d-c +2$, then
$$abcd \leq (a+1)(b-1)(c-1)(d+1),$$
where the inequality is strict iff $b-a > d-c+2$.
\end{lem}
\begin{proof}
We have 
\begin{gather}\label{rel2}
(a+1)(b-1)(c-1)(d+1) = (ab + b - a - 1)(cd + c  - d -1).
\end{gather}
and by assumption
\begin{gather}\label{rel1}
b-a -1 \geq d - c +1 
\end{gather}
This gives
$$
\begin{array}{ccl}
(a+1)(b-1)(c-1)(d+1) - abcd  & = & cd(b-a-1) - ( ab +  b- a -1)(d-c+1)\\
 \text{ because of } \eqref{rel1}  &\geq & cd(b-a-1) - ( ab +  b- a -1)(b-a-1)\\
\text{ because of } a < c, b < d & \geq &((a+1)(b+1)- ab - (b-a-1))(b-a - 1)  \\
& = & (2a +2)(b-a-1)\geq 0\\
\end{array}
$$
since $b-a \geq 1$.
\end{proof}

\begin{thm}\label{C2} Let $\mathcal{P}$ be of rank $2$ and $\lambda \in \mathcal{P}^+$. Further let $\lie g$ be a simple complex classical Lie algebra of rank $2$, e.g. $\lie g$ is of type $A_2$ or $C_2$. If $\blambda = (\lambda_1, \lambda_2) \trianglelefteq \bmu  = (\mu_1, \mu_2) \in \mathcal{P}^+(\lambda,2)$, then
$$\dim V(\iota(\lambda_1)) \otimes V(\iota(\lambda_2)) \leq \dim V(\iota(\mu_1)) \otimes V(\iota(\mu_2)).$$
\end{thm}

\begin{proof}
The map $\iota$ is an isomorphism for $\lie g$ of type $A_2$ or $C_2$, hence we will identify $\lambda$ with $\iota(\lambda)$ throughout the proof to simplify the notation.\\

\noindent
Let $\lie g$ be of type $A_2$, then 
\[\blambda \trianglelefteq \bmu \Leftrightarrow \blambda \preceq \iota \bmu \]
and the theorem follows from \cite[Theorem 1 (i)]{CFS2012}.\\
Let $\lie g$ be of type $C_2$, then we have the set of positive coroots $\{h_1, h_2, h_{1,2}, h_{1, \overline{2}} \}$ (see Remark~\ref{coroot-list}).\\
It suffices to show the theorem for $\bmu = (\mu_1, \mu_2)$ being a possible cover of $\blambda$. Let $\sigma \in S_n$ such that $\sigma(\lambda_1 - \lambda_2) \in \mathcal{P}^+$, by Proposition~\ref{cover-elements} we may assume that $\sigma \in \{ id, \bos_{1,2}, \bos_{2,3}\}$.
Write $\lambda_1$ as $a \omega_1 + b \omega_2$ then  $\lambda_2 = (n-a) \omega_1 + (m - b) \omega_2$, where $n = \lambda_1(h_1), m = \lambda_1(h_2)$.
We will use Weyl's dimension formula and recall that $\iota^{-1}(\rho) = \omega_1 + \omega_2 \in \mathcal{P}^+$. It suffices to  show the following inequalities: 
\begin{eqnarray}
 \langle \lambda_1, h_1 \rangle \,  \langle \lambda_2, h_1 \rangle \, \langle \lambda_1, h_{1, \overline{2}} \rangle \,  \langle \lambda_2, h_{1, \overline{2} } \rangle \, \leq \langle \mu_1, h_1 \rangle \,  \langle \mu_2, h_1 \rangle \, \langle \mu_1, h_{1, \overline{2}} \rangle \,  \langle \mu_2, h_{1, \overline{2} } \rangle 
\label{eq:C2-1}
\end{eqnarray}
and for $h_\alpha \in \{ h_2, h_{1,2} \}$
\begin{eqnarray}
\langle \lambda_1, h_\alpha \rangle \,  \langle \lambda_2, h_\alpha \rangle \leq \langle \mu_1 , h_\alpha \rangle \, \langle \mu_2 , h_\alpha \rangle. \label{eq:C2-2}
\end{eqnarray}

\noindent
Both inequalities together give
\[ \prod_{\alpha \in R^+} \langle \lambda_1, h_\alpha \rangle \, \langle \lambda_2, h_\alpha \rangle \leq \prod_{\alpha \in R^+} \langle \mu_1, h_\alpha  \rangle \, \langle \mu_2, h_\alpha \rangle \]
this finishes then by (\ref{eq:dim-tensor}) the proof.

\noindent
In \cite[Proof of Theorem 1(i)]{CFS2012} it was shown that for a given coroot $h_\alpha$
\begin{eqnarray}
\bor_{\alpha, \ell} (\iota(\blambda)) & \leq & \bor_{\alpha, \ell} (\iota(\bmu)) \text{ for all } 1\leq \ell \leq 2
\end{eqnarray}
\begin{eqnarray}
\Rightarrow \; \; \langle \lambda_1, h_\alpha \rangle \,  \langle \lambda_2, h_\alpha \rangle \leq \langle \mu_1 , h_\alpha \rangle \, \langle \mu_2 , h_\alpha \rangle.
\label{AN-case}
\end{eqnarray}

\noindent
This gives \eqref{eq:C2-2} and the rest of the proof is dedicated to proving \eqref{eq:C2-1}. We will show this by case by case considerations for $\sigma \in \{\id, \bos_{1,2} \}$, $\bos_{2,3}$ is similar. Throughout the proof we will use the following simple fact
\begin{eqnarray}
x < y \in \bz_{>  0} \Rightarrow (x+1)(y-1) \leq xy \text{ with equality iff } x = y-1.
\label{eq:simple-fact}
\end{eqnarray}

\noindent
We replace in the following $\lambda_i$ by $\lambda_i + \rho$, $\mu_i$ by $\mu_i + \rho$ to avoid having a summand $\rho(h_\alpha)$ in each factor. This does not change the argument, the possible covers for this shifted tuple are exactly the same. Note that 
\[\lambda_1  (h_1) = a, \lambda_1 (h_{\alpha}) = a + 2b, \lambda_2(h_1) = n-a, \lambda_2(h_{\alpha}) = n-a + 2(m-b).\]

\noindent
If $\sigma = \id$, then $a\geq n-a, b \geq m-b$, the three possible covers (two of type I and one of type II, Corollary~\ref{number-cover}) are  (see Proposition~\ref{cover-elements})
\begin{enumerate}
\item $(\mu_1, \mu_2) = (\lambda_1 - \omega_1, \lambda_2+ \omega_1)$ (type I). \\
		By assumption 
		\[n-a < a \text{ and }n-a+1 + 2(m-b) < n-a + 2(m-b),\] hence with (\ref{eq:simple-fact}) we have $(a-1)(n-a+1) \geq a(n-a)$ and 
		\[(a-1 + 2b)(n-a+1 + 2(m-b)) \geq (a + 2b)(n-a + 2(b-m)).\] 
                 Combining both we have as desired 
		\[a(n-a)(a+2b)(n-a + 2(m-b))\]
		\[ \leq (a-1)(n-a+1)(a-1 + 2b)(n-a+1 + 2(m-b)).\]
\item $(\mu_1, \mu_2) = (\lambda_1 - \omega_2, \lambda_2+ \omega_2)$ is similar.\\
\item $ (\mu_1, \mu_2) = (a \omega_1 + (m-b) \omega_2, (n-a) \omega_1 + b \omega_2)$ (type II).\\ 
		By assumption 
		\[a + 2b \geq n-a + 2b \geq (n-a) + 2(m-b)\]
		and
		\[  a + 2b \geq a + 2(m-b)  \geq (n-a) + 2(m-b).\]
		With (\ref{eq:simple-fact}) it follows
		\[(a+2b)(n-a + 2(m-b)) \leq (a + 2(m-b))(n-a + 2b)\]
		this then implies as desired
		\[[a(n-a)(a + 2b)][(n-a + 2(m-b))] \; 
		\leq \; [a(n-a)] [(a + 2(m-b)) (n-a + 2b)] . \]
\end{enumerate}

\noindent
If $\sigma = \bos_{1,2}$, then $a < n-a, 2b - m \geq n -2a > 0$, the three possible covers are
\begin{enumerate}
     \item  $(\mu_1, \mu_2) = (\lambda_1 - \bos_{1,2}(\omega_1), \lambda_2 + \bos_{1,2}( \omega_1))$ (type I).\\
		Here we must distinguish two cases
       		\begin{enumerate}
   		     \item Suppose $a + 2b \geq n-a +2(m-b)$, then (see (\ref{eq:simple-fact})
				\[(a +1 +2(b -1))(n-a -1  + 2(m-b+1)) \geq (a +2b)(n-a + 2(m-b))\] 
				and by assumption 
				\[ (a+1)(n-a-1) \geq a(n-a).\]
				Combining both proves the claim.
      		     \item Suppose $a + 2b \leq n-a + 2(m-b)$. By assumption we have $a < n-a$, so to apply Lemma~\ref{4-tuple} it suffices to show:
				\[ (n-2a) - ((n-2a) + 2(m-2b)) = - 2(m-2b) \geq  2.\] 
				But this follows, since by assumption $2b - m \geq 1$, so the claim follows
    	        \end{enumerate}
     \item  $(\mu_1, \mu_2) = (\lambda_1 - \omega_2, \lambda_2 + \omega_2)$ (type I).\\
		 Then $a + 2b \geq n-a + 2(m-b) +2$ (since $2b > m$ and $a < n-a$) and so by (\ref{eq:simple-fact})
		\[(a + 2b)((n-a + 2(m-b)))\leq (a + 2(b-1)) ( n-a+  2 (m- b+1)),\] 
		this implies the claim.
     \item $(\mu_1, \mu_2) = (\lambda_1 - \bos_1( (n-2a) (\omega_1)), \lambda + \bos_1((n-2a) (\omega_1)))$ (type II).\\ 
		Then 
		\[a + 2b \geq n-a + 2(b + 2a -n) \geq (n-a + 2(m-b)),\]
		where the first inequality is due to $a < n-a$ and the second inequality due to $2b - m < n -2a$. Further
		\[a + 2b \geq a + 2(m+n-b-2a) \geq (n-a + 2(m-b)),\]
		where the first inequality is due to $2b -m \geq n- 2a$ and the second inequality due to $a < n-a$. Then (\ref{eq:simple-fact}) implies
     		\[(a+2b) (n-a + 2(m-b)) \leq (n-a + 2(b + 2a - n)) (a + 2(m+n-b-2a)).\]
\end{enumerate}

\noindent
We omit the similar computations for $\sigma = \bos_{2,3}$.     
\end{proof}

\noindent
Let $\lie g$ be of type $B_3$ or $D_4$, recall the embedding $\iota: \mathcal{P}^+ \longrightarrow P^+$, so the weights we are considering are supported on the first two nodes only. 
\begin{cor}\label{bd-cor}
Let $\blambda = (\lambda_1, \lambda_2) \trianglelefteq \bmu  = (\mu_1, \mu_2) \in \mathcal{P}^+(\lambda,2)$, then
$$\dim V(\iota(\lambda_1)) \otimes V(\iota(\lambda_2)) \leq \dim V(\iota(\mu_1)) \otimes V(\iota(\mu_2)).$$
\end{cor}
\begin{proof}
Let $\lie g$ be of type $B_3$ and $\lambda \in \mathcal{P}^+$, then $\iota(\lambda)$ is supported on the first two nodes only, hence Weyl's dimension formula reduces to $\dim V(\iota(\lambda)) = $
\[ ( \langle \iota(\lambda), h_2 \rangle^2 \, \langle \iota(\lambda), h_1 + h_2 \rangle^2 \, 
\langle 2 \iota(\lambda) , h_2 \rangle \, \langle 2 \iota(\lambda), h_1 + h_2 \rangle \, \langle \iota(\lambda), h_1 \rangle \, \langle \iota(\lambda), h_1 + 2 h_2 \rangle )\prod_{\alpha \in R^+} \frac{1}{(\rho, h_\alpha)}.
 \]
We want to show again the inequality 
\[ \dim V(\iota(\lambda_1)) \otimes V(\iota(\lambda_2)) \leq \dim V(\iota(\mu_1)) \otimes V(\iota(\mu_2)).\]
From \eqref{AN-case} we know that under the assumption $\blambda \trianglelefteq \bmu$
\[\langle \iota(\lambda_1), h \rangle \, \langle \iota(\lambda_2), h \rangle \leq  \langle\iota(\mu_1), h \rangle \, \langle \iota(\mu_2), h \rangle \]
for all $h \in \{ h_1, h_2, h_{1,2} =  h_1 + h_2 \}$. So it remains to prove the inequality for $h_{1,\overline{2}} = h_1 + 2 h_2$. But this follows from \eqref{eq:C2-1}.

\noindent
Let $\lie g$ be of type $D_4$, then Weyl's dimension formula reduces to
$$
\dim V(\iota(\lambda)) = ( \langle \iota(\lambda), h_2 \rangle^4 \, \langle \iota(\lambda), h_1 + h_2 \rangle^4 \,  \langle \iota(\lambda), h_1 \rangle \, \langle \iota(\lambda), h_1 + 2 h_2 \rangle ) \prod_{\alpha \in R^+} \frac{1}{(\rho, h_\alpha)} .
$$
and the analogous argument gives the proof here.
\end{proof}

\subsection{}
Here we will finally prove Theorem~\ref{mainthm}, for this let $\lambda \in \mathcal{P}^+$ and $\lie g$ be of type $A_n, C_n, B_{n+1}, D_{n+2}$. We want to show that if $\blambda \trianglelefteq \bmu \in \mathcal{P}^+(\lambda,2)$, then
$$\dim V(\iota(\lambda_1)) \otimes V(\iota(\lambda_2)) \leq \dim V(\iota(\mu_1)) \otimes V(\iota(\mu_2)).$$
For this we will use again Weyl's dimension formula
$$\dim V(\iota(\lambda)) = \prod_{\alpha \in R^+} \frac{\langle \iota(\lambda) , h_\alpha \rangle}{(\rho, h_\alpha)}.$$
and show a general form of the inequalities (\ref{eq:C2-1}) and  (\ref{eq:C2-2}), this will prove the claim of the theorem.

\noindent
If $\lie g$ is of type $A_n$, then this follows immediately from \cite[Theorem 1 (i)]{CFS2012}.\\

\noindent
If $\lie g$ is of type $C_n$ and $h_\alpha$ a coroot such that there exists $i \in I$ with $\omega_i(h_\alpha)= 2$, we then say $h_\alpha$ has height $2$. Then (see Remark~\ref{coroot-list})
$$h_\alpha = h_i + \ldots + h_{j-1} + 2 (h_j + \ldots + h_n).$$
for some $i < j$. The Lie algebra $\lie g'$ associated to the coroots
\[h_{i,j-1} = h_i + \ldots + h_{j-1}, h_{j,n} = h_j + \ldots + h_n\] 
(that is the Lie algebra generated by the corresponding root vectors) is of type $C_2$, the set of the corresponding positive coroots is $ \{h_{i,j-1}, h_{j,n}, h_{i, j-1} + h_{j,n}, h_{i, j-1} + 2 h_{j,n} \}$.\\
Since $\blambda \trianglelefteq \bmu$, we have 
\[\bor_{(i,j-1),1}(\blambda) \leq \bor_{(i,j-1)}(\bmu), r_{(j,n), 1}(\blambda) \leq \bor_{(j,n), 1}(\bmu).\] 
Then the same arguments as in the proof of Theorem~\ref{C2} gives (see \ref{eq:C2-1}) that
\[\langle \iota(\lambda_1),  h_{i,j-1} \rangle \,  \langle \iota(\lambda_1),  h_{\alpha} \rangle \langle \iota(\lambda_2),  h_{i,j-1} \rangle \langle \iota(\lambda_2),  h_{i,j-1}+ 2 h_{j,n}\rangle \]
\[ \leq  \langle \iota(\mu_1),  h_{i,j-1} \rangle \langle \iota(\mu_1),  h_{\alpha} \rangle \langle \iota(\mu_2),  h_{i,j-1} \rangle \langle \iota(\mu_2),  h_{i,j-1} + 2 h_{j,n}\rangle. \]

\noindent
If we subtract from the set of all coroots the set $\{ h_{i, j-1}, h_{i,j-1} + 2 h_{j,n} \, | \, i  < j \in I \}$ we are left with coroots of height $1$.\\
It remains to show that for all such $h_\alpha$ we have
$$\langle \iota(\lambda_1), h_\alpha \rangle \langle \iota(\lambda_2), h_\alpha \rangle \leq \langle \iota(\mu_1), h_\alpha \rangle \langle \iota(\mu_2), h_\alpha \rangle.$$
This follows as \ref{eq:C2-2} from the proof of \cite[Theorem 1(i)]{CFS2012}.\\

\noindent
If $\lie g$ is of type $B_{n+1}$ and $h_\alpha$ is a root of height $2$, then 
$$h_\alpha = h_{i,j-1} + 2 h_{j,n} + h_{n+1},$$
where $1 \leq i \leq j \leq n$, and we set $h_{j,j-1} = 0$. 
First of all, if $i = j-1$, then $\iota(\lambda)(h_\alpha) = 2 \iota(\lambda)(h_{j,n})$ (since $\iota(\lambda)$ is not supported on $h_{n+1}$). Since  $\blambda \trianglelefteq \bmu$ implies that $r_{\alpha,1}(\iota(\blambda)) \leq r_{\alpha,1}(\iota(\bmu))$, the proof of \cite[Theorem 1(i)]{CFS2012} gives the inequality for this coroot.\\

\noindent
If $i < j-1$ then the Lie algebra generated by the root vectors corresponding to $h_{i,j-1}, h_{j,n}, h_{n+1}$ is a simple Lie algebra of type $B_3$. As in the $C_2$ case above: $\blambda \trianglelefteq \bmu$ implies
\[\bor_{(i,j-1),1}(\blambda) \leq \bor_{(i,j-1)}(\bmu), r_{(j,n), 1}(\blambda) \leq \bor_{(j,n), 1}(\bmu).\] 
Then Corollary~\ref{bd-cor} gives as for (\ref{eq:C2-1})
\[\langle \iota(\lambda_1), h_\alpha \rangle \langle \iota(\lambda_1), h_{i,j-1} \rangle \langle \iota(\lambda_2), h_\alpha \rangle \langle \iota(\lambda_2), h_{i,j-1} \rangle \]
\[\leq
\langle \iota(\mu_1), h_\alpha \rangle \langle \iota(\mu_1), h_{i,j-1} \rangle \langle \iota(\mu_2), h_\alpha \rangle \langle \iota(\mu_2), h_{i,j-1}.\]
The remaining coroots have all height $1$, so the needed inequality follows again as in \cite[Theorem 1(i)]{CFS2012}.

\noindent
If $\lie g$ is of type $D_{n+2}$, and $h_\alpha$ of height $2$. Then as in the $C_2$ (resp. $B_3$) case, we have an induced simple Lie algebra of type $D_4$ and again Corollary~\ref{bd-cor} gives as for (\ref{eq:C2-1}) the needed inequality here. Again the remaining coroots have all height $1$ and the inequalities follow again as in  \cite[Theorem 1(i)]{CFS2012}.

\noindent
All together this implies in the several cases
$$\prod_{\alpha \in R^+} \langle \iota(\lambda_1), h_\alpha \rangle \langle \iota(\lambda_2), h_\alpha \rangle \leq \prod_{\alpha \in R^+} \langle \iota(\mu_1), h_\alpha \rangle \langle \iota(\mu_2), h_\alpha \rangle$$
this gives the proof of Theorem~\ref{mainthm}.

\subsection{} To conclude the paper, it remains to prove (\ref{eq:max-dim}), so if  $\bmu = (\mu_1, \ldots, \mu_k) \in \mathcal{P}^+(\lambda, k)$  we have
\begin{eqnarray}
\dim V(\iota(\mu_1)) \otimes \ldots \otimes V(\iota(\mu_k)) \leq \dim V(\iota(\lambda_1)) \otimes \ldots \otimes V(\iota(\lambda_k)) 
\end{eqnarray}
where $\blambda^{\max} = (\lambda_1,\ldots, \lambda_k)$ and equality iff $\bmu \sim \blambda^{\max}$, hence iff $\bmu$ lies in the $S_k$-orbit of $\blambda^{\max}$.

\noindent
Suppose $\bmu$ is not in the $S_k$-orbit of $\blambda^{\max}$, then (\ref{eq:max-uni-2}) implies that there exists $i < j$ such that $(\mu_i, \mu_j)$ is not maximal in $\mathcal{P}^+(\mu_i + \mu_j, 2)$.  Denote the maximal $2$-tuple in $\mathcal{P}^+(\mu_i + \mu_j,2)$ by $(\nu_i, \nu_j)$. 
Then we have by Theorem~\ref{mainthm}
\[ \dim V(\iota(\mu_i)) \otimes V(\iota(\mu_j)) <  \dim V(\iota(\nu_i)) \otimes V(\iota(\nu_j)).\]
We define a new $k$-tuple 
\[\bnu = (\nu_1, \ldots, \nu_k)\]
where we set  $\nu_{\ell} = \mu_{\ell}$ for $\ell \neq i,j$. Then Proposition~\ref{prop:extend-tuple} implies $\bmu \triangleleft \bnu$ and we have
\[ \dim V(\iota(\mu_1)) \otimes \ldots \otimes V(\iota(\mu_k)) < \dim V(\iota(\nu_1)) \otimes \ldots \otimes V(\iota(\nu_k)) .\]
Again by (\ref{eq:max-uni-2}) we have that $\blambda^{\max}$ is uniquely determined by the condition that $(\lambda_i, \lambda_j)$ is maximal in $\mathcal{P}^+(\lambda_i+ \lambda_j, 2)$ for all $i < j$. So there is an increasing chain of tuples, such that in each step only two components are changed, so each step is of the form $\bmu \triangleleft \bnu$ as above. An induction along this chain gives now (\ref{eq:max-dim}).

\end{document}